# Some Comparisons for Gaussian Processes[*]

Richard A. Vitale
Department of Statistics, U-120
University of Connecticut
Storrs, CT 06269–3120
rvitale@uconnvm.uconn.edu.

**Abstract**

Extensions and variants are given for the well-known comparison principle of Gaussian processes based on ordering by pairwise distance.

## 1 Introduction

Among the most important tools for Gaussian processes are comparison principles, or simply *comparisons*. Typically they provide for the derivation of upper bounds by majorizing a given process with a second process that is larger in some sense as well as having more tractable properties. For general discussions, see Adler (1990), Fernique (1997), Ledoux and Talagrand (1991), and Lifshits (1995). The purpose of this note is to elaborate some variants of perhaps the most widely applied comparison:

**Theorem 1** *Suppose that* $\{X_i, i \in I\}$ *and* $\{Y_i, i \in I\}$ *are two mean–zero Gaussian processes indexed by the same denumerable set* $I$, *and supposea that*

$$E(X_i - X_j)^2 \leq E(Y_i - Y_j)^2 \tag{1}$$

*for all* $i, j \in I$. *Then for any non–decreasing, convex* $g : I\!R_+ \to I\!R^1$,

$$Eg\left(\sup_{i,j}(X_i - X_j)\right) \leq Eg\left(\sup_{i,j}(Y_i - Y_j)\right) \tag{2}$$

*and*

$$E \sup_i X_i \leq E \sup_i Y_i. \tag{3}$$

This is associated with Sudakov (1971, 1976) and Fernique (1975) with a later proof by Alexander (1985) (Ledoux and Talagrand also mention unpublished work of S. Chevet). Important extensions and variants have been given by Gordon (1985, 1987, 1992) and Kahane (1986).

We wish to show a sequence of extensions that are apparently new. **In the sequel, $\{X_i, i \in I\}$ and $\{Y_i, i \in I\}$ continue to be mean–zero Gaussian processes indexed by a**

---
[*]*Proc. Amer. Math. Soc.* **128** (2000), 3043–3046.

denumerable $I$ and satisfying (1); $\{m_i, i \in I\}$ are arbitrary constants. The first extension is

$$E \sup_i \{X_i + m_i\} \leq E \sup_i \{Y_i + m_i\}. \tag{4}$$

This implies that, for any constant $m$ and $k \in I$,

$$E \sup_i \{X_i + m_i - X_k, m\} \leq E \sup_i \{Y_i + m_i - Y_k, m\}, \tag{5}$$

Finally, *for any non–decreasing, convex* $g : \mathbb{R}^1 \to \mathbb{R}^1$ *with (i)* $g(-\infty) > -\infty$ *or with (ii)* $\max\left[Eg_+\left(\sup_i \{X_i + m_i\} - X_k\right), Eg_+\left(\sup_i \{Y_i + m_i\} - Y_k\right)\right] < \infty$, *one has*

$$Eg\left(\sup_i \{X_i + m_i\} - X_k\right) \leq Eg\left(\sup_i \{Y_i + m_i\} - Y_k\right). \tag{6}$$

In the next section, we provide proofs. The third section shows how a special case is related to a set of weakened Slepian–Schläfli assumptions.

## 2 Proofs

We use the fact that any positive constant can be regarded as nearly the supremum of a centered Gaussian process (*cf.* Vitale, 1996). A technical formulation and proof are as follows:

**Lemma 1** *For any* $c \geq 0$, *there is a sequence* $\{W_j\}_{j=1}^\infty$ *of mean zero Gaussian variables such that almost surely (i)* $c \leq \sup_{j \geq n} W_j < \infty$ *and (ii)* $\sup_{j \geq n} W_j \downarrow c$ *as* $n \to \infty$.

**Proof** By homogeneity, it is enough to consider $c = 1$. Referring to an example of Marcus and Shepp (1972), one has the following: for independent, standard Gaussian variables $\{Z_j\}_1^\infty$ and $\{\sigma_j\}_1^\infty$ with $1/\sigma_j^2 = 2 \log j + 2 \log \log j$, let $W_1 = \sigma_1 Z_1, W_2 = -W_1, W_3 = \sigma_2 Z_2, W_4 = -W_3, \cdots$. Then $P\left(\sup_{j \geq 1} W_j \geq 1\right) = 1$ and $\sup_{j \geq 1} W_j$ has an atom at 1. In fact, the required verifications for these two properties show that they depend only on the tail behavior of the sequence $\{\sigma_j\}$. Hence, more generally, $P\left(\sup_{j \geq n} W_j \geq 1\right) = 1$ and $\sup_{j \geq n} W_j$ has an atom at 1. A standard 0-1 argument using the tail measurability of $\limsup_j W_j$ then yields the required convergence.

**Proofs of (4), (5), (6).** To establish (4), it suffices to consider $I$ finite, say $I = \{1, 2, \ldots, N\}$, since the general case follows from the Monotone Convergence Theorem; then one may also assume that $0 \leq \min_i m_i$ since adding a constant to each side of (4) amounts to a shifted set of $m_i$. Finally, without loss of generality and using the lemma, one may assume that the underlying probability space is rich enough to support mean-zero Gaussian variables $\{W_{ij}\}$ that are independent of the $\{X_i, Y_i\}$, and such that for each $1 \leq i \leq N$, $\sup_{j \geq n} W_{ij} \geq m_i$ and $\sup_{j \geq n} W_{ij} \downarrow m_i$ as $n \to \infty$.

Consider

$$\sup_i \left\{ X_i + \sup_{j \geq n} W_{ij} \right\} = \sup_{i, j \geq n} \{X_i + W_{ij}\}.$$

A similar expression for the $Y$ process can also be written as as a supremum of mean zero Gaussian variables. It can be checked that the two collections of variables are ordered according to (1). Consequently, (3) implies

$$E \sup_i \left\{ X_i + \sup_{j \geq n} W_{ij} \right\} \leq E \sup_i \left\{ Y_i + \sup_{j \geq n} W_{ij} \right\}. \tag{7}$$

Each integrand is bounded from below by 0 and is decreasing in $n$, so that letting $n \to \infty$ and applying the Dominated Convergence Theorem then yields (4).

For (5), set $\hat{X}_i = X_i - X_k$ and $\hat{Y}_i = Y_i - Y_k$ for $i \neq k$ with $\hat{X}_k = \hat{Y}_k = 0$; further set $\hat{m}_i = m_i - m_k$ for $i \neq k$ and $\hat{m}_k = \max\{m, m_k\}$. Then (4) holds for the hatted system, and this is equivalent to (5).

For (6), assume as before that $I$ is finite and $g(-\infty) > -\infty$. From (5), the following holds for $-\infty < t < \infty$:

$$E \sup_i \{X_i + m_i - X_k, t\} \leq E \sup_i \{Y_i + m_i - Y_k, t\}.$$

Subtracting $t$ leads to the equivalent form

$$E g_t \left( \sup_i \{X_i + m_i - X_k\} \right) \leq E g_t \left( \sup_i \{Y_i + m_i - Y_k\} \right),$$

where $g_t : \mathbb{R}^1 \to \mathbb{R}^1$ is given by $g_t(\cdot) = (\cdot - t)_+$. Then it is enough to recall that the closed, positive linear hull of the collection $\{g_t\}$ is precisely the class of non–negative, non–decreasing, convex functions on $\mathbb{R}^1$; adding a constant adjusts for a given limit at $-\infty$.

For the alternate assumption, note that the truncation $\max\{g, c\}$ conforms to the previous requirement and that the truncation can be removed by letting $c \to -\infty$ and appealing to the Monotone Convergence Theorem.

## 3   A Connection with the Slepian–Schläfli Comparison

The well–known Slepian–Schläfli comparison ([12], [13]) gives a stronger statement than (3) at the expense of a more stringent hypothesis. It can be formulated as follows:

*Suppose that, together with (1), one has for all $i \in I$*

$$E X_i^2 = E Y_i^2. \tag{8}$$

*for all $t > 0$. Then*

$$P(\sup_i X_i > t) \leq P(\sup_i Y_i > t). \tag{9}$$

Unfortunately, the requirement of strict equality in (8) often precludes applicability. While there seems to be no easy remedy for this, let us show that our earlier results can be adapted to a weak form of (8), which leads in turn to an integrated form of (9):

$$\int_t^\infty P\left( \sup_i X_i > s \right) ds \leq \int_t^\infty P\left( \sup_i Y_i > s \right) ds \tag{10}$$

for all $t > 0$. In place of (8), consider

$$EX_i^2 \leq EY_i^2 \tag{11}$$

for all $i \in I$. Then, assuming (1) and (11), we have

$$E \sup_i \{X_i + m_i, m\} \leq E \sup_i \{Y_i + m_i, m\}. \tag{12}$$

and *for any non–decreasing, convex $g : \mathbb{R}^1 \to \mathbb{R}^1$ with (i) $g(-\infty) > -\infty$ or with (ii)* $\max \{Eg_+ (\sup_i \{X_i + m_i\}), Eg_+ (\sup_i \{Y_i + m_i\})\} < \infty$,

$$Eg \left( \sup_i \{X_i + m_i\} \right) \leq Eg \left( \sup_i \{Y_i + m_i\} \right). \tag{13}$$

These can be argued as follows, where it suffices once again to assume that $I = \{1, 2, \cdots, N\}$. Consider the augmented process $\{X_1, X_2, \cdots, X_N\} \cup \{X_0 = 0\}$. Having assumed (1) and (11) for the original process corresponds precisely to (1) holding for the augmented process. We then choose $k = 0$, and let $m_0 = -\infty$ so that $X_i + m_i|_{i=0} = m_0$ and $Y_i + m_i|_{i=0} = m_0$ do not participate in any suprema. Then (5) and (6) translate to (12) and (13), the analogue of (4) reducing to a special case of (13).

Finally, (10) follows from (13):

$$\int_t^\infty P \left( \sup_i X_i > s \right) ds = E \left( \sup_i X_i - t \right)_+ \leq E \left( \sup_i Y_i - t \right)_+ = \int_t^\infty P \left( \sup_i Y_i > s \right) ds.$$

## 4 Acknowledgment

The referee and Probability Editor provided helpful comments.